\documentclass[preprint,3p,sort&compress,final,times]{elsarticle}
\usepackage{amsmath}
\usepackage{amssymb}
\usepackage{mathtools}
\usepackage{placeins}
\usepackage[usenames,dvipsnames]{color}
\usepackage{xcolor}
\usepackage{cases}
\usepackage{array}
\usepackage{tabto}
\usepackage{algorithm}
\usepackage[noend]{algpseudocode}

\usepackage[usenames,dvipsnames]{color}

\journal{}
\begin{document}
\begin{frontmatter}

\title{High order tracer variance stable transport with low order energy conserving dynamics for the thermal shallow water equations}
\author[BOM]{David Lee\corref{cor}}
\ead{david.lee@bom.gov.au}
\author[ANU]{Kieran Ricardo}
\author[MON]{Tamara Tambyah}
\address[BOM]{Bureau of Meteorology, Melbourne, Australia}
\address[ANU]{Mathematical Sciences Institute, Australian National University, Canberra, Australia}
\address[MON]{School of Mathematics, Monash University, Melbourne, Australia}
\cortext[cor]{Corresponding author. Tel. +61 452 262 804.}

\begin{abstract}
	A high order discontinuous Galerkin method for the material transport of thermodynamic tracers is coupled
	to a low order mixed finite element solver in the context of the thermal shallow water equations. The coupling
	preserves the energy conserving structure of the low order dynamics solver, while the high order material
	transport scheme is provably tracer variance conserving, or damping with the inclusion of upwinding.
	The two methods are coupled via a nested hierarchy of meshes, with the low order mesh of the dynamics 
	solver being embedded within the high order transport mesh, for which the basis functions are collocated at 
	the Gauss-Legendre quadrature points.

	Standard test cases are presented to verify the consistency and conservation properties of the method. While the overall
	scheme is limited by the formal order of accuracy of the low order dynamics, the use of high order, tracer variance conserving
	transport is shown to preserve richer turbulent solutions without compromising model stability compared to a purely low order
	method.
\end{abstract}

\end{frontmatter}

Low order numerical schemes are common for implicit atmospheric solvers using finite difference \cite{Wood14,Yeh02}, finite volume \cite{Sandbach15}
and finite element \cite{Melvin19} methods, since the lower condition number of the corresponding Jacobian operator results in faster convergence
and matrix assembly compared to high order methods. The approximate Jacobian operators used in these methods typically linearise around the 
fast acoustic, buoyancy and gravity modes \cite{Maynard20,Lee24}, which are not resolved in time, 
and so there is less emphasis on the accuracy of their representation compared to the slow moving vortical and inertial motions. 
Conversely transport terms that explicitly resolve these inertial motions are typically based on high order methods
for the construction of 
numerical fluxes using finite volume \cite{SG11,BK25} or discontinuous Galerkin \cite{Bendall19} methods, since the 
upwinding of low order methods as required in order to suppress high frequency artefacts leads to excessive dissipation at 
low order.

In the present work we harmonise these two approaches within a consistent, energy and tracer variance conserving formulation using
a mixed finite element method for the low order dynamics and a high order discontinuous Galerkin method for the material transport
of thermodynamic tracers. We do this within the context of the thermal shallow water equations, where the buoyancy (the thermodynamic
variable) is represented in flux form for the low order dynamics. The low order buoyancy fluxes are computed with respect to the high order
buoyancy as derived via discontinuous Galerkin material transport. This same high order representation of the buoyancy is also
used in the low order pressure gradient term in the momentum equation, so as to preserve the anti-symmetry of the flux form buoyancy transport
and the pressure gradient, and thus satisfy energy conservation \cite{Eldred19}. Meanwhile the high order discontinuous Galerkin representation
of the material transport of buoyancy is expressed as a combination of both the weak and strong form advection operators so as to
provably conserve tracer variance, and hence improve model stability \cite{Ricardo24a,Ricardo24b,Tambyah25,Ricardo25}, since tracer variance 
constitutes a numerical entropy of the dynamical system. However in contrast to these previous tracer variance conserving formulations,
here we express the tracer variance conserving discrete advection operator in material form and not flux form, so as to couple this with 
the low order flux form dynamics.

The integration of the high order discontinuous Galerkin thermodynamic transport into the low order dynamics solver is achieved by 
constructing the Lagrange polynomial basis functions so as to be orthogonal with respect to the interior Gauss-Legendre quadrature 
points \cite{KG10,GK11}, while the low order mesh is constructed such that the cell centres are coincident with these high order quadrature
points, such that there is a 1:1 mapping between the high and low order degrees of freedom. For the low order discretisation, this
is done with a hierarchy of coarser meshes, such that a single
high order discontinuous Galerkin element, with four degrees of freedom in each dimension, is coincident with a single low order element
two levels down in the mesh refinement hierarchy from the one on which the dynamics are represented. 

The remainder of this article proceeds as follow: In Section 1 we introduce the thermal shallow water equations, which are used as
a model of a geophysical system with a distinction between wave-like and inertial processes for the remainder of this article, and 
their conservation properties. In Section 2 we present the low order dynamics, high order advection variational discretisation 
with reference to the conservation properties described in Section 1. In Section 3 we present results for standard test cases
to verify the method with respect to its consistency and conservation properties, and in Section 4 we discuss the conclusions drawn 
from this study.

\section{Energy and tracer variance conservation for the thermal shallow water equations}

The rotating thermal shallow water equations may be expressed for the velocity, $\boldsymbol{u}$, depth $h$ and depth 
weighted buoyancy, $S=sh$ (with $s$ being the buoyancy) over the two-dimensional domain $\Omega\subset\mathbb{R}^2$ as 
\cite{Eldred19,KLZ21}
\begin{subequations}\label{eq::tsw1}
	\begin{align}
		\frac{\partial\boldsymbol{u}}{\partial t} + (\nabla\times\boldsymbol{u} + f)\cdot\boldsymbol{u}^{\perp} + 
		\frac{1}{2}\nabla(\boldsymbol{u}\cdot\boldsymbol{u} + sh) + \frac{1}{2}s\nabla h &=0,\label{eq::tsw1::u}\\
		\frac{\partial h}{\partial t} + \nabla\cdot(h\boldsymbol{u}) &= 0,\label{eq::tsw1::h}\\
		\frac{\partial S}{\partial t} + \nabla\cdot(sh\boldsymbol{u}) &= 0,\label{eq::tsw1::S}
	\end{align}
\end{subequations}
where $f$ is the Coriolis term and $\boldsymbol{u}^{\perp} = (-v,u)$ for the vector velocity $\boldsymbol{u}=(u,v)$. 
Alternatively, the buoyancy transport may be formulated as the material transport of $s$, instead of the flux form transport of $S$. 
Expanding \eqref{eq::tsw1::S} and invoking \eqref{eq::tsw1::h} we have
\begin{equation}
	h\frac{\partial s}{\partial t} + s\frac{\partial h}{\partial t} + h\boldsymbol{u}\cdot\nabla s + s\nabla\cdot(h\boldsymbol{u}) = 
	h\Bigg(\frac{\partial s}{\partial t} + \boldsymbol{u}\cdot\nabla s\Bigg) = 0,
\end{equation}
such that
\begin{equation}
	\frac{\partial s}{\partial t} + \boldsymbol{u}\cdot\nabla s = 0
\end{equation}
may instead be used in place of \eqref{eq::tsw1::S}.

Using $S$ as the prognostic variable for the buoyancy, the corresponding energy is defined over the domain $\Omega\subset\mathbb{R}^2$ as
\begin{equation}
	\mathcal{H} = \frac{1}{2}\int h\boldsymbol{u}\cdot\boldsymbol{u} + hS\mathrm{d}\Omega,
\end{equation}
for which the variational derivatives are given with respect to the prognostic variables, $\boldsymbol{u},h,S$, as:
\begin{subequations}
	\begin{align}
		\frac{\delta\mathcal{H}}{\delta\boldsymbol{u}} &= h\boldsymbol{u} := \boldsymbol{F},\\
		\frac{\delta\mathcal{H}}{\delta h} &= \frac{1}{2}(\boldsymbol{u}\cdot\boldsymbol{u} + S) := \Phi,\\
		\frac{\delta\mathcal{H}}{\delta S} &= \frac{h}{2} := T.
	\end{align}
\end{subequations}
The original system \eqref{eq::tsw1} may then be re-formulated in skew-symmetric non-canonical Hamiltonian form with respect to the 
variational derivatives as
\begin{subequations}\label{eq::tsw2}
	\begin{align}
		\frac{\partial\boldsymbol{u}}{\partial t} + q\boldsymbol{F}^{\perp} + 
		\nabla\Phi + s\nabla T &=0,\label{eq::tsw2::u}\\
		\frac{\partial h}{\partial t} + \nabla\cdot\boldsymbol{F} &= 0,\label{eq::tsw2::h}\\
		\frac{\partial S}{\partial t} + \nabla\cdot(s\boldsymbol{F}) &= 0,
	\end{align}
\end{subequations}
where $q = (\nabla\times\boldsymbol{u} + f)/h$ is the potential vorticity. Energy conservation is established by left multiplication
of \eqref{eq::tsw2} by the variational derivatives and integration by parts such that (for periodic or wall boundary conditions)
\begin{multline}\label{eq::en_cons_cont}
	\int\frac{\mathrm{d}\mathcal{H}}{\mathrm{d}t}\mathrm{d}\Omega = 
	\int\frac{\delta\mathcal{H}}{\delta\boldsymbol{u}}\cdot\frac{\partial\boldsymbol{u}}{\partial t} + 
	\frac{\delta\mathcal{H}}{\delta h}\frac{\partial h}{\partial t} + 
	\frac{\delta\mathcal{H}}{\delta S}\frac{\partial S}{\partial t}\mathrm{d}\Omega = 
	\int-\boldsymbol{F}\cdot(q\boldsymbol{F}^{\perp} + \nabla\Phi + s\nabla T) 
	- \Phi\nabla\cdot\boldsymbol{F} - T\nabla\cdot(s\boldsymbol{F})\mathrm{d}\Omega = 0.
\end{multline}

In addition to the energy, \eqref{eq::tsw1} also conserves the tracer variance, 
\begin{equation}
	\mathcal{S} = \frac{1}{2}\int\frac{S^2}{h}\mathrm{d}\Omega,
\end{equation}
for which the variational derivatives are given as
\begin{subequations}
	\begin{align}
		\frac{\delta\mathcal{S}}{\delta\boldsymbol{u}} &= \boldsymbol{0},\\
		\frac{\delta\mathcal{S}}{\delta h} &= -\frac{1}{2}\frac{S^2}{h^2} = -\frac{1}{2}s^2,\\
		\frac{\delta\mathcal{S}}{\delta S} &= \frac{S}{h} = s.
	\end{align}
\end{subequations}
Left multiplication of \eqref{eq::tsw2} by the variational derivatives of $\mathcal{S}$ as given above gives as an analogue
to the energy conservation in \eqref{eq::en_cons_cont}
\begin{multline}\label{eq::tv_cons_cont}
	\int\frac{\mathrm{d}\mathcal{S}}{\mathrm{d}t}\mathrm{d}\Omega = 
	\int\frac{\delta\mathcal{S}}{\delta\boldsymbol{u}}\cdot\frac{\partial\boldsymbol{u}}{\partial t} + 
	\frac{\delta\mathcal{S}}{\delta h}\frac{\partial h}{\partial t} + 
	\frac{\delta\mathcal{S}}{\delta S}\frac{\partial S}{\partial t}\mathrm{d}\Omega = 
	\int-\boldsymbol{0}\cdot(q\boldsymbol{F}^{\perp} + \nabla\Phi + s\nabla T) 
	+ \frac{1}{2}s^2\nabla\cdot\boldsymbol{F} - s\nabla\cdot(s\boldsymbol{F})\mathrm{d}\Omega = \\
	\int-\frac{1}{2}s^2\nabla\cdot\boldsymbol{F} - \frac{1}{2}\boldsymbol{F}\cdot\nabla s^2\mathrm{d}\Omega
	=0.
\end{multline}
For the purpose of deriving a discrete formulation that conserves a discrete analogue of $\mathcal{S}$, it is important
to note that while energy conservation, as given in \eqref{eq::en_cons_cont} is dependent only on integration by
parts in space and the chain rule in time (subject to appropriate boundary conditions), tracer variance conservation
as given in \eqref{eq::tv_cons_cont} is also dependent on the product rule. Since the product rule is challenging 
to satisfy discretely in the absence of continuous basis functions and exact integration \cite{LPG18}, the discrete 
tracer variance conserving formulation derived in the proceeding Section 2.2 will be
tailored specifically to negate the need for this by first expanding the flux terms at the continuous level.

\section{Discrete formulation}

\subsection{Mixed finite element formulation of the thermal shallow water equations}

Mixed compatible variational formulations for the thermal shallow water equations (which preserve integration by parts 
discretely and hence conserve energy) have been presented previously \cite{Eldred19,Ricardo24a,Tambyah25}, so these
are discussed only briefly here. We introduce the discrete function spaces of the form $\mathbb{W}_0^L\subset H^1(\Omega)$,
$\mathbb{W}_1^L\subset H(\mathrm{div},\Omega)$, $\mathbb{W}_2^L\subset L^2(\Omega)$, which consist of polynomial functions 
that are square integrable over the domain $\Omega$ with respect to the $H^1$, $H(\mathrm{div})$ and $L^2$ norms respectively, 
and restricted to the lowest polynomial order of these spaces. In practice for $\Omega\subset\mathbb{R}^2$ this means polynomials
that are piecewise linear and $C^0$ continuous in both dimensions for $\mathbb{W}_0^L$, vector polynomials that are
piecewise constant and discontinuous in the tangent direction and piecewise linear and $C^0$ continuous in the normal direction
for $\mathbb{W}_1^L$ and discontinuous and piecewise constant for $\mathbb{W}_2^L$.

Introducing the test functions $\boldsymbol{v}_h\in\mathbb{W}_1^L$, $\phi_h,\sigma_h\in\mathbb{W}_2^L$ and integrating
over the domain $\Omega$ and the time step $\Delta t$, we seek solutions for the discrete analogues of the prognostic variables 
$\boldsymbol{u}_h\in\mathbb{W}_1^L$, $h_h,S_h\in\mathbb{W}_2^L$ via the discrete form of \eqref{eq::tsw2} as
\begin{subequations}\label{eq::tsw_disc_1}
\begin{align}
		\int\boldsymbol{v}_h\cdot(\boldsymbol{u}_h^{n+1} - \boldsymbol{u}_h^n)\mathrm{d}\Omega + 
		\Delta t\int\boldsymbol{v}_h\cdot\overline{q}_h\overline{\boldsymbol{F}}_h^{\perp}\mathrm{d}\Omega \notag\\ 
		-\Delta t\int\nabla\cdot\boldsymbol{v}_h\overline{\Phi}_h\mathrm{d}\Omega + 
		\Delta t\int\boldsymbol{v}_h\cdot\overline{s}_h\nabla \overline{T}_h\mathrm{d}\Omega - 
		\Delta t\int\boldsymbol{v}_h\cdot\hat{\boldsymbol{n}}\{\overline{s}_h\}[\overline{T}_h]\mathrm{d}\Gamma &= 0,\label{eq::tsw_disc::u}\\
		\int\phi_h(h_h^{n+1} - h_h^n)\mathrm{d}\Omega + 
		\Delta t\int\phi_h\nabla\cdot\overline{\boldsymbol{F}}_h\mathrm{d}\Omega &= 0,\label{eq::tsw_disc::h}\\
		\int\sigma_h(S_h^{n+1} - S_h^n)\mathrm{d}\Omega - 
		\Delta t\int\nabla\sigma_h\cdot \overline{s}_h\overline{\boldsymbol{F}}_h\mathrm{d}\Omega + 
		\Delta t\int[\sigma_h]\{\overline{s}_h\}\overline{\boldsymbol{F}}_h\cdot\hat{\boldsymbol{n}}\mathrm{d}\Gamma &= 0
		\label{eq::tsw_disc::S},
\end{align}
\end{subequations}
where $\Gamma\subset\mathbb{R}$ denotes the edge facets of the elements, and $\hat{\boldsymbol{n}}$ represents the outward unit 
normal at the element boundary. The operators $[a]:=a^+-a^-$, $\{a\}:=(a^++a^-)/2$ represent the jump and mean operators respectively
(with $a^+$ and $a^-$ being derived from the cells in the positive and negative direction with respect to the orientation of 
$\hat{\boldsymbol{n}}$ respectively). The $\overline{a}$ operator denotes exact temporal integration of $a$ over a discrete time level 
$\Delta t = t^{n+1} - t^n$.

Note that we have applied integration by parts to the discrete analogue of the Bernoulli potential term, $\nabla\Phi$ in \eqref{eq::tsw2::u}. 
Since $\boldsymbol{v}_h\in\mathbb{W}_1^2$ is div-conforming, no boundary 
integral term arises from this, whereas the weak form of the discrete analogue of the term $\nabla\cdot(s\boldsymbol{F})$ does generate
a boundary integral since $s\boldsymbol{F}$ is not div-conforming. Also note that while this term has been presented in the weak form,
the discrete form of the term $s\nabla T$ has been presented in the strong form \cite{KG10} so as to discretely preserve the anti-symmetry
and thus the energy conservation of the pressure gradient and buoyancy flux terms. Since 
$\boldsymbol{v}_h,\boldsymbol{F}_h\in\mathbb{W}_1^L$ are $C^0$ continuous in the direction of $\hat{\boldsymbol{n}}$, their normal 
components are uniquely defined along $\Gamma$.

For the low order discretisation, functions in $\mathbb{W}_2^L$ are piecewise constant, and hence their derivatives vanish. Therefore 
the terms $\Delta t\int\boldsymbol{v}_h\cdot\overline{s}_h\nabla \overline{T}_h\mathrm{d}\Omega$ and 
$\Delta t\int\nabla\sigma_h\cdot \overline{s}_h\overline{\boldsymbol{F}}_h\mathrm{d}\Omega$ may be omitted from 
\eqref{eq::tsw_disc::u} and \eqref{eq::tsw_disc::S}, since these vanish for $\overline{T}_h,\sigma_h\in\mathbb{W}_2^L$.

The discrete variational derivatives and potential vorticity are computed exactly to second order in time
\cite{Eldred19,BC18,LP21,Tambyah25} between time levels $n$ and $n+1$ for $\overline{\boldsymbol{F}}_h\in\mathbb{W}_1^L$,
$\overline{\Phi}_h,\overline{T}_h\in\mathbb{W}_2^L$, $\overline{q}_h\in\mathbb{W}_0^L$ as
\begin{subequations}\label{eq::tsw_diag_1}
	\begin{align}
		\int\boldsymbol{v}\cdot\overline{\boldsymbol{F}}_h\mathrm{d}\Omega &= 
		\frac{1}{6}\int\boldsymbol{v}_h\cdot(2h_h^n\boldsymbol{u}_h^{n} + h_h^n\boldsymbol{u}_h^{n+1} + 
		h_h^{n+1}\boldsymbol{u}_h^{n} + 2h_h^{n+1}\boldsymbol{u}_h^{n+1})\mathrm{d}\Omega, \\
		\int\phi_h\overline{\Phi}_h\mathrm{d}\Omega &= 
		\int\phi_h\Bigg(\frac{1}{6}(\boldsymbol{u}_h^{n}\cdot\boldsymbol{u}_h^{n} + 
		\boldsymbol{u}_h^{n}\cdot\boldsymbol{u}_h^{n+1} + \boldsymbol{u}_h^{n+1}\cdot\boldsymbol{u}_h^{n+1}) 
		+ \frac{1}{4}(S_h^n + S_h^{n+1})\Bigg)
		\mathrm{d}\Omega, \\ 
		\int\sigma_h\overline{T}_h\mathrm{d}\Omega &= \frac{1}{4}\int\sigma_h(h_h^n + h_h^{n+1})\mathrm{d}\Omega, \\
		\frac{1}{2}\int\psi_h(h_h^n+h_h^{n+1})\overline{q}_h\mathrm{d}\Omega &= 
		\frac{1}{2}\int\nabla^{\perp}\psi_h\cdot(\boldsymbol{u}_h^n + \boldsymbol{u}_h^{n+1})\mathrm{d}\Omega + 
		\int\psi_h f\mathrm{d}\Omega,
	\end{align}
\end{subequations}
where $\psi_h\in\mathbb{W}_0^L$. The mean value of the low order buoyancy over the time level, $\overline{s}_h\in\mathbb{W}_2^L$
will be discussed in Section 2.3.

Energy is conserved discretely in space and time subject to the exact integration over the time level for the variational derivatives 
$\overline{\boldsymbol{F}}_h,\overline{\Phi}_h,\overline{T}_h$ by assigning these to $\boldsymbol{v}_h,\phi_h,\sigma_h$ respectively
in \eqref{eq::tsw_disc_1} and summing all terms in a discrete analogue of \eqref{eq::en_cons_cont} \cite{CH11,Eldred19,Tambyah25}.
In addition to the energy, \eqref{eq::tsw_disc_1} also conserves the total mass, $\int h_h\mathrm{d}\Omega$ and total depth weighted
buoyancy, $\int S_h\mathrm{d}\Omega$, which are are satisfied by setting $\phi_h=1$ and 
$\sigma_h=1$ in \eqref{eq::tsw_disc::h} and \eqref{eq::tsw_disc::S} respectively. As for the discrete energy conservation, discrete
density weighted buoyancy conservation is satisfied for any choice of $\overline{s}_h$, so this is preserved for the 
high order buoyancy transport discussed below.

\subsection{Discrete tracer variance conservation in material form}

In order to derive a discrete tracer variance material transport expression for the buoyancy, we begin by first expanding the flux form
expression at the continuous level \eqref{eq::tsw1::S} as
\begin{equation}\label{eq::adv_S_cont}
	\frac{\partial S}{\partial t} + \frac{1}{2}\Bigg(\nabla\cdot(\boldsymbol{F} s) + 
	\boldsymbol{F}\cdot\nabla s +  s\nabla\cdot\boldsymbol{F}\Bigg) = 0.
\end{equation}
This continuous form is the starting point for previous discrete tracer variance conserving formulations for
flux form transport \cite{Ricardo24a,Ricardo24b,Tambyah25}. To derive an analogous expression for the 
material form transport of $s$, we expand the time derivative for $S=hs$ in \eqref{eq::adv_S_cont} and recall 
the continuity equation \eqref{eq::tsw2::h}, such that
\begin{equation}
	h\frac{\partial s}{\partial t} +  s\Bigg(\frac{\partial h}{\partial t} + 
	\nabla\cdot\boldsymbol{F}\Bigg) + \frac{1}{2}\Bigg(\nabla\cdot(\boldsymbol{F} s) + 
	\boldsymbol{F}\cdot\nabla s -  s\nabla\cdot\boldsymbol{F}\Bigg) = 0.
\end{equation}
In order to derive a discrete, tracer variance conserving analogue to the material transport expression above, 
we begin by introducing the high order space $\mathbb{W}_2^H\subset L^2(\Omega)$ which is spanned by the set of square integrable
polynomials of degree $p$ over the domain $\Omega$ and discontinuous along the element boundaries $\Gamma$. 
Introducing the high order test function $\chi_h\in\mathbb{W}_2^H$, eliminating 
the continuity equation (which at the discrete level is satisfied pointwise as in \eqref{eq::tsw_disc::h}) and integrating over 
the domain $\Omega$ with respect to $\chi_h$ gives
\begin{equation}
	\int\chi_h h_h\frac{\partial s_h}{\partial t}\mathrm{d}\Omega + 
	\frac{1}{2}\int\chi_h\Bigg(\nabla\cdot(\boldsymbol{F}_h s_h) + 
	\boldsymbol{F}_h\cdot\nabla s_h -  s_h\nabla\cdot\boldsymbol{F}_h\Bigg)\mathrm{d}\Omega = 0.
\end{equation}
While $\boldsymbol{F}_h\cdot\hat{\boldsymbol{n}}$ is continuous and div-conforming over the element boundaries, 
$(\boldsymbol{F}_h s_h)\cdot\hat{\boldsymbol{n}}$ and $(\nabla s_h)\cdot\hat{\boldsymbol{n}}$ are not. Hence we 
integrate by parts for the former term and apply boundary 
integrals to penalise against discontinuities in both terms as
\begin{multline}\label{eq::buoy_adv_disc}
	\int\chi_h h_h\frac{\partial s_h}{\partial t}\mathrm{d}\Omega
	- \frac{1}{2}\int\nabla\chi_h\cdot\boldsymbol{F}_h s_h\mathrm{d}\Omega
	+ \frac{1}{2}\int [\chi_h]\{s_h\}\boldsymbol{F}_h\cdot\hat{\boldsymbol{n}}\mathrm{d}\Gamma \\
	+ \frac{1}{2}\int\chi_h\boldsymbol{F}_h\cdot\nabla s_h\mathrm{d}\Omega
	- \frac{1}{2}\int \boldsymbol{F}_h\cdot\hat{\boldsymbol{n}}\{\chi_h\}[ s_h]\mathrm{d}\Gamma
	- \frac{1}{2}\int\chi_h s_h\nabla\cdot\boldsymbol{F}_h\mathrm{d}\Omega 
	+ \alpha\int|\boldsymbol{F}_h\cdot\hat{\boldsymbol{n}}|[\chi_h][ s_h]\mathrm{d}\Gamma
	= 0.
\end{multline}
Note that we have also added the upwinding term $\int|\boldsymbol{F}_h\cdot\hat{\boldsymbol{n}}|[\chi_h][ s_h]\mathrm{d}\Gamma$
which is enabled for the choice of the upwinding parameter as $\alpha=1$.

For a pointwise divergence free mass flux, $\nabla\cdot\boldsymbol{F}_h=0$, and constant in time fluid depth, $\partial h_h/\partial t=0$, 
\eqref{eq::buoy_adv_disc} will discretely conserve 
total buoyancy such that $\frac{\mathrm{d}}{\mathrm{d}t}\int h_hs_h\mathrm{d}\Omega = \int h_h\frac{\partial s_h}{\partial t}\mathrm{d}\Omega = 0$. 
This is assured for the choice of test 
function as $\chi_h=1$, which can be represented exactly in the discrete space of $\mathbb{W}_2^H$, giving
\begin{multline}\label{eq::adv_mass_con}
	\int h_h\frac{\partial s_h}{\partial t}\mathrm{d}\Omega
	- \frac{1}{2}\int\nabla 1\cdot\boldsymbol{F}_h s_h\mathrm{d}\Omega
	+ \frac{1}{2}\int [1]\{s_h\}\boldsymbol{F}_h\cdot\hat{\boldsymbol{n}}\mathrm{d}\Gamma 
	+ \frac{1}{2}\int \boldsymbol{F}_h\cdot\nabla s_h\mathrm{d}\Omega
	- \frac{1}{2}\int \boldsymbol{F}_h\cdot\hat{\boldsymbol{n}}\{1\}[ s_h]\mathrm{d}\Gamma \\
	- \frac{1}{2}\int s_h\nabla\cdot\boldsymbol{F}_h\mathrm{d}\Omega 
	+ \alpha\int|\boldsymbol{F}_h\cdot\hat{\boldsymbol{n}}|[1][ s_h]\mathrm{d}\Gamma
		= \\
	\int h_h\frac{\partial s_h}{\partial t}\mathrm{d}\Omega
	+ \frac{1}{2}\int \boldsymbol{F}_h\cdot\nabla s_h\mathrm{d}\Omega
	- \frac{1}{2}\int \boldsymbol{F}_h\cdot\hat{\boldsymbol{n}}[ s_h]\mathrm{d}\Gamma
	- \frac{1}{2}\int s_h\nabla\cdot\boldsymbol{F}_h\mathrm{d}\Omega 
		= \\
	\int h_h\frac{\partial s_h}{\partial t}\mathrm{d}\Omega
		- \int s_h\nabla\cdot\boldsymbol{F}_h\mathrm{d}\Omega = 
		\int h_h\frac{\partial s_h}{\partial t}\mathrm{d}\Omega = 0.
\end{multline}
In addition, \eqref{eq::buoy_adv_disc} also discretely conserves the second moment, the tracer variance, which is given as
\begin{equation}\label{eq::S}
	\mathcal{S}_h = \frac{1}{2}\int h_h s_h^2\mathrm{d}\Omega,
\end{equation}
for which the variational derivatives are given as
\begin{align}
	\int\phi_h\frac{\delta\mathcal{S}_h}{\delta h_h}\mathrm{d}\Omega &= \frac{1}{2}\int\phi_h s_h^2\mathrm{d}\Omega, \\
	\int\chi_h\frac{\delta\mathcal{S}_h}{\delta s_h}\mathrm{d}\Omega &= \int\chi_h h_h s_h\mathrm{d}\Omega,
\end{align}
where $\phi_h$ is the low order test function in the same space as $h_h$. Since \eqref{eq::tsw_disc::h} holds pointwise we are not 
restricted to test functions in $\mathbb{W}_2^L$, and are free to chose higher order test functions for this.
Setting the test functions for the low order continuity equation 
and the high order transport equations respectively as $\phi_h=\frac{1}{2} s_h^2$, $\chi_h= s_h$ and adding 
the two expressions gives
\begin{multline}
	\frac{1}{2}\int s_h^2\Bigg(\frac{\partial h}{\partial t} + \nabla\cdot\boldsymbol{F}_h\Bigg)\mathrm{d}\Omega +
	\int s_h h_h\frac{\partial s_h}{\partial t}\mathrm{d}\Omega
	- \frac{1}{2}\int\nabla s_h\cdot\boldsymbol{F}_h s_h\mathrm{d}\Omega
	+ \frac{1}{2}\int [ s_h]\{s_h\}\boldsymbol{F}_h\cdot\hat{\boldsymbol{n}}\mathrm{d}\Gamma \\
	+ \frac{1}{2}\int s_h\boldsymbol{F}_h\cdot\nabla s_h\mathrm{d}\Omega
	- \frac{1}{2}\int \boldsymbol{F}_h\cdot\hat{\boldsymbol{n}}\{s_h\}[ s_h]\mathrm{d}\Gamma
	- \frac{1}{2}\int s_h^2\nabla\cdot\boldsymbol{F}_h\mathrm{d}\Omega 
	+ \alpha\int|\boldsymbol{F}_h\cdot\hat{\boldsymbol{n}}|[ s_h][ s_h]\mathrm{d}\Gamma
	= 0.
\end{multline}
Cancellation and integration by parts (assuming continuous time) yields
	\begin{equation}\label{eq::disc_var_con}
	\frac{\mathrm{d}\mathcal{S}_h}{\mathrm{d} t} = 
	\frac{1}{2}\int\frac{\partial( h_h s_h^2)}{\partial t}\mathrm{d}\Omega = 
	-\alpha\int|\boldsymbol{F}_h\cdot\hat{\boldsymbol{n}}|[ s_h]^2\mathrm{d}\Gamma \le 0.
\end{equation}
such that the tracer variance is conserved for a choice of $\alpha=0$ (a centered flux) and dissipated for $\alpha>0$. 
For $\alpha=1$ the flux is fully upwinded (taking contributions from the upstream cell only), whereas for $0<\alpha<1$
the flux is only partially upwinded.

Designing solvers to ensure the conservation of additional high order invariants (beyond the energy) in time
for non-canonical Hamiltonian systems is challenging, as it is difficult 
to preserve the temporal chain rule discretely as in \eqref{eq::en_cons_cont}. Here we use a stiffly stable 
third order Runge-Kutta integrator (SSP-RK3) \cite{SO88} (which will not conserve tracer variance in time) 
for the temporal integration of \eqref{eq::buoy_adv_disc},
with the low order depth and mass fluxes derived from their time averages, $\overline{h}_h$, 
$\overline{\boldsymbol{F}}_h$. While it is mathematically desirable to conserve tracer variance in both space
and time, in practice the absence of some form of tracer variance dissipation such as the upwinding term presented
in \eqref{eq::buoy_adv_disc} leads to excessive grid-scale noise due to non-linear aliasing and a loss of
coherence for fine scale features \cite{Ricardo24a,Tambyah25}.

Since the depth and mass flux are derived from the low order dynamics,
the overall scheme is limited to second order accuracy. However as will be shown in Section 3.3, the
use of \eqref{eq::buoy_adv_disc} leads to stable solutions with less dissipation of tracer variance
than would be the case if $s_h$ was derived purely from the low order solution. 

\subsection{Coupling of low order dynamics and high order advection schemes}

As detailed above, the conservation of tracer variance for high order methods on discontinuous function 
spaces requires the skew-symmetric expansion of the advection operator at the continuous level since the product 
rule is not satisfied discretely \cite{Ricardo24a,Ricardo24b,Tambyah25}. However for low order methods using
a piecewise constant representation of the tracer $\overline{s}_h\in\mathbb{W}_2^L$ \eqref{eq::tsw_disc_1}, 
the tracer gradient vanishes within the element and only the boundary integrals remain, such that spatial conservation
of tracer variance is satisfied directly. The low order tracer variance for the flux form variables is given as
$\mathcal{S}_h^L(h_h,S_h) = \frac{1}{2}\int S_h^2/h_h\mathrm{d}\Omega$.
Setting the test functions as the variational derivatives of $\mathcal{S}_h^L$ such that
$\phi_h = -\frac{1}{2}\overline{s}^2_h$, $\sigma_h=\overline{s}_h$ in 
\eqref{eq::tsw_disc::h}, \eqref{eq::tsw_disc::S} respectively for $\overline{s}_h\in\mathbb{W}_2^L$ and adding gives
\begin{subequations}\label{eq::lo_tv_disc}
	\begin{align}
	-\frac{1}{2}\int\overline{s}_h^2\frac{\partial h_h}{\partial t}\mathrm{d}\Omega + 
		\int\overline{s}_h\frac{\partial S_h}{\partial t}\mathrm{d}\Omega &= 
	\frac{1}{2}\int\overline{s}_h^2\nabla\cdot\overline{\boldsymbol{F}}_h\mathrm{d}\Omega + 
	\int\nabla\overline{s}_h\cdot\overline{s}_h\overline{\boldsymbol{F}}_h\mathrm{d}\Omega -
	\int\overline{\boldsymbol{F}}_h\cdot\hat{\boldsymbol{n}}[\overline{s}_h]\{\overline{s}_h\}\mathrm{d}\Gamma \\
		&= 
	\int\overline{\boldsymbol{F}}\cdot\hat{\boldsymbol{n}}[\overline{s}_h]\{\overline{s}_h\}\mathrm{d}\Gamma -
		\int\overline{\boldsymbol{F}}\cdot\hat{\boldsymbol{n}}[\overline{s}_h]\{\overline{s}_h\}\mathrm{d}\Gamma \label{eq::lo_tv_disc::2} \\
		&= 0,
	\end{align}
\end{subequations}
where in the first term of \eqref{eq::lo_tv_disc::2} we recall the identity
$[\overline{s}_h^2] = 2\{\overline{s}_h\}[\overline{s}_h]$ and integration by parts.

The low order mixed finite element dynamics solver for $\boldsymbol{u}_h^{n+1},h_h^{n+1},S_h^{n+1}$ 
\eqref{eq::tsw_disc_1}, \eqref{eq::tsw_diag_1},
and the high order discontinuous Galerkin advection solver for $s_h^{n+1}$ \eqref{eq::buoy_adv_disc} are coupled
through 
a hierarchy of sub-divided meshes, with the high order mesh for the buoyancy transport being two levels of 
refinement higher than the coarse level mesh (four degrees of freedom in each dimension instead of one).
At the first level of refinement, the cell centers are shifted such that these are collocated with the Gauss-Legendre
quadrature points at the next level of refinement. This ensures that the element centres for each ``patch'' of 
$4\times 4$ low order elements is collocated with the Gauss-Legendre points of a single high order element of polynomial degree $p=3$.

At each time step $n$ the initial condition for the high order buoyancy $s_h^{n}$ is 
derived from the low order dynamics for $\gamma_h\in\mathbb{W}_2^H$ as
\begin{equation}\label{eq::s_diag}
	\int\gamma_h h_h^ns_h^n\mathrm{d}\Omega = \int\gamma_h S_h^n\mathrm{d}\Omega.
\end{equation}
The high order buoyancy at the new time level $n+1$ at nonlinear Newton iteration $k$, $s_h^{n+1}$ is then 
derived from the high order discontinuous Galerkin scheme \eqref{eq::buoy_adv_disc}
using the low order mass flux and depth fields averaged over the time interval,
$\overline{\boldsymbol{F}}_h,\overline{h}_h$.

The time averaged low order buoyancy in \eqref{eq::tsw_disc::u}, \eqref{eq::tsw_disc::S} is then given 
as $\overline{s}_h=\frac{1}{2}(\Pi(s_h^n)+\Pi(s_h^{n+1}))$, 
where $\Pi(\cdot):\mathbb{W}_2^H\rightarrow\mathbb{W}_2^L$ represents a projection from the high order 
discontinuous Galerkin space in which $s_h^{n+1}$ is computed to the low order space in which $\overline{s}_h$
is represented. For the four point Lagrange polynomials collocated with the Gauss-Legendre quadrature points,
and the low order fine scale cells centered at those same points, $\Pi(\cdot)$ is simply a diagonal matrix
consisting of the Jacobian determinant within each fine scale low order element.

In order to negate the 
inertial CFL limit of the explicit buoyancy transport, one could also sub-step the buoyancy transport,
as is done in some operational models \cite{Melvin19,BK25}. 

The mixed problem is solved using two iterations of GMRES \cite{SS86} 
using the PETSc library \cite{petsc-web-page,petsc-user-ref,petsc-efficient}. This is not an
efficient strategy for large domains in parallel as unlike geometric multigrid methods, which use Jacobi \cite{Sandbach15,Maynard20} 
or patch based smoother approaches \cite{Arnold00,Hiptmair97}, GMRES
and other Krylov methods require global norm and dot products that do not scale efficiently with the 
number of processors. 
For the mixed velocity-pressure problem Jacobi iteration is less effective
owing to the more complex structure of the eigenvalues, and so 
more sophisticated smoothers are required that account for kernels in the vector field space 
\cite{Arnold00}, such as decomposition of the vector field into its rotational and divergent components
\cite{Hiptmair97} or hybridised methods \cite{Betteridge23} in the context of mixed finite elements,
and patch based additive Schwarz schemes in the context of collocated finite volume methods \cite{YC14}.

The quasi-Newton problem for which the low order dynamics solver is applied is given as
\begin{equation}\label{eq::newton}
	\begin{bmatrix}
		\boldsymbol{\mathsf{M}}_1 & -\frac{\Delta t}{4}g\boldsymbol{\mathsf{D}}^{\top} & 
		-\frac{\Delta t}{4}\boldsymbol{\mathsf{D}}^{\top} \\ 
		\frac{\Delta t}{2}H\boldsymbol{\mathsf{D}} & \boldsymbol{\mathsf{M}}_2 & \boldsymbol{\mathsf{0}} \\
		\boldsymbol{\mathsf{0}} & \boldsymbol{\mathsf{0}} & \boldsymbol{\mathsf{M}}_2 
	\end{bmatrix}
	\begin{bmatrix}
		\delta\boldsymbol{u}_h^k\\ \delta h_h^k\\ \delta S_h^k
	\end{bmatrix} = 
	-\begin{bmatrix}
		\mathcal{R}_{\boldsymbol{u}}\\ \mathcal{R}_h\\ \mathcal{R}_S
	\end{bmatrix},
\end{equation}
where the residual vectors $\mathcal{R}_{\boldsymbol{u}},\mathcal{R}_h,\mathcal{R}_S$ represent the prognostic 
equations \eqref{eq::tsw_disc_1}, subject to the diagnostic terms \eqref{eq::tsw_diag_1} at each Newton iteration,
and $g,H$ are the mean values of the buoyancy and depth.
The current estimate of the solution for time level $n+1$ at the end of each Newton iteration $k$ is then updated as 
$\boldsymbol{u}_h^{n+1} = \boldsymbol{u}_h^{n+1} + \delta\boldsymbol{u}_h^k$, $h_h^{n+1}=h_h^{n+1} + \delta h_h^k$, 
$S_h^{n+1} = S_h^{n+1} + \delta S_h^k$, until the system is converged to below some specified tolerance. The matrix 
operators in \eqref{eq::newton} are given as
\begin{subequations}
	\begin{align}
		\boldsymbol{\mathsf{M}}_1 &:= \int\boldsymbol{v}_h\cdot\boldsymbol{w}_h\mathrm{d}\Omega
		\qquad\forall\boldsymbol{v}_h,\boldsymbol{w}_h\in\mathbb{W}_1^L,\\
		\boldsymbol{\mathsf{M}}_2 &:= \int\phi_h\sigma_h\mathrm{d}\Omega\qquad\forall\phi_h,\sigma_h\in\mathbb{W}_2^L,\\
		\boldsymbol{\mathsf{D}} &:= \int\phi_h\nabla\cdot\boldsymbol{v}_h\mathrm{d}\Omega
		\qquad\forall\boldsymbol{v}_h\in\mathbb{W}_1^L,\phi_h\in\mathbb{W}_2^L.
	\end{align}
\end{subequations}
The approximate Jacobian in \eqref{eq::newton} is similar to that used previously for a mixed variational
form of the thermal shallow water equations with material buoyancy transport \cite{Eldred19}, and a simplified
version of that previously used for flux form buoyancy transport \cite{Tambyah25}. To improve convergence
one may replace the mean values of buoyancy and depth, $g$ and $H$ with spatially varying values from some
previous time step or Newton iteration.
The full time stepping scheme can then be expressed in Algorithm 1, where $\epsilon$ represents the tolerance
and $k_{max}$ the maximum number of nonlinear iterations as prescribed for the solver.
\makeatletter
\def\BState{\State\hskip-\ALG@thistlm}
\makeatother
\begin{algorithm}
\caption{Coupled low order dynamics high order transport time stepping}\label{euclid}
\begin{algorithmic}[1]
	\For{time step $n$}
		\State$\mathrm{diagnose}:\ s_h^n\in\mathbb{W}_2^H\ \mathrm{via}\ \eqref{eq::s_diag}$
		\State$\mathrm{set}:\ \boldsymbol{u}_h^{n+1}=\boldsymbol{u}_h^n,\ h_h^{n+1}=h_h^n,\ S_h^{n+1} = S_h^n$
		\For{Newton iteration $k$}
		\State$\mathrm{diagnose}:\ \overline{\boldsymbol{F}}_h,\overline{\Phi}_h,\overline{T}_h,\overline{q}_h\ \mathrm{via}\ \eqref{eq::tsw_diag_1}$
		\State$\mathrm{integrate}:\ s_h^{n+1}\in\mathbb{W}_2^H\ \mathrm{via}\ \eqref{eq::buoy_adv_disc}$
		\State$\mathrm{diagnose}:\ \overline{s}_h = \frac{1}{2}(\Pi(s_h^n)+\Pi(s_h^{n+1}))\in\mathbb{W}_2^L$
		\State$\mathrm{solve\ for}:\ \delta\boldsymbol{u}_h^k,\ \delta h_h^k,\ \delta S_h^k\ \mathrm{via}\ \eqref{eq::tsw_disc_1}\ \mathrm{and}\ \eqref{eq::newton}$
		\State$\mathrm{update }:\ \boldsymbol{u}_h^{n+1} = \boldsymbol{u}_h^{n+1} + \delta\boldsymbol{u}_h^k,\ h_h^{n+1} = h_h^{n+1} + \delta h_h^k,\ S_h^{n+1} = S_h^{n+1} + \delta S_h^k$
		\If{$|\delta\boldsymbol{u}_h^k|/|\boldsymbol{u}_h^{n+1}| < \epsilon\ \mathrm{and}\ |\delta h_h^k|/|h_h^{n+1}| < \epsilon\ \mathrm{and}\ |\delta S_h^k|/|S_h^{n+1}| < \epsilon\ \mathrm{\mathbf{or}}\ k = k_{max}$} $\mathrm{break}$
		\EndIf
		\EndFor
	\EndFor
\end{algorithmic}
\end{algorithm}

\section{Results}

\subsection{Advection only: solid body rotation}

In order to verify the high order discontinuous Galerkin scheme in \eqref{eq::buoy_adv_disc} and 
its conservation properties, we first apply this to a stand alone test case of solid body rotation 
within a periodic domain of size $\Omega=[-\pi,\pi)\times[-\pi,\pi)$ with a constant analytical
depth of $h=1$ and a constant analytical mass flux field of $\boldsymbol{F} = (y,-x)$.
The tracer field is initialised with a state of 
$s_h^0 = \exp(-10((x-x_o)^2 + (y-y_o)^2))$, where $(x_o,y_o)=(-0.4\pi,0.4\pi)$ is the initial position of the tracer,
and periodic boundary conditions are applied.
The tracer makes a single revolution in $2\pi$ time units, after which the $L^2$ errors are computed with 
respect to its initial position. The domain is discretised using $12,24,48$ and $96$ elements using cubic 
polynomials (degree $p=3$), and time steps of $\Delta t=\pi/300,\pi/600,\pi/1200$ and $\pi/2400$. Two different 
configurations are presented, one with upwinded fluxes, using 
$\alpha=1$ in \eqref{eq::buoy_adv_disc}, and a second using centered fluxes ($\alpha=0$).

\begin{figure}[!hbtp]
\begin{center}
\includegraphics[width=0.48\textwidth,height=0.36\textwidth]{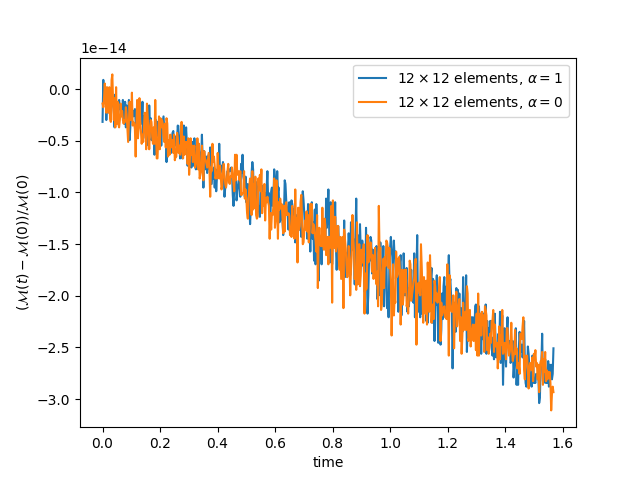}
\includegraphics[width=0.48\textwidth,height=0.36\textwidth]{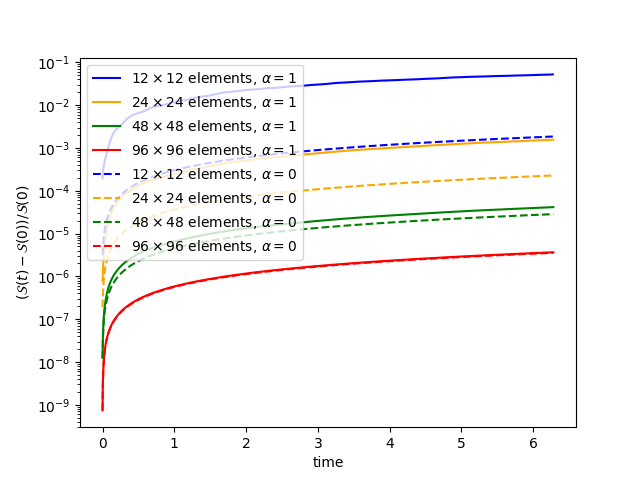}
	\caption{Mass (left) and tracer variance (right) conservation errors with time for the high order discontinuous Galerkin 
	material advection scheme using upwinded ($\alpha=1$) and centered ($\alpha=0$) fluxes.}
\label{fig::adv1}
\end{center}
\end{figure}

\begin{figure}[!hbtp]
\begin{center}
\includegraphics[width=0.48\textwidth,height=0.36\textwidth]{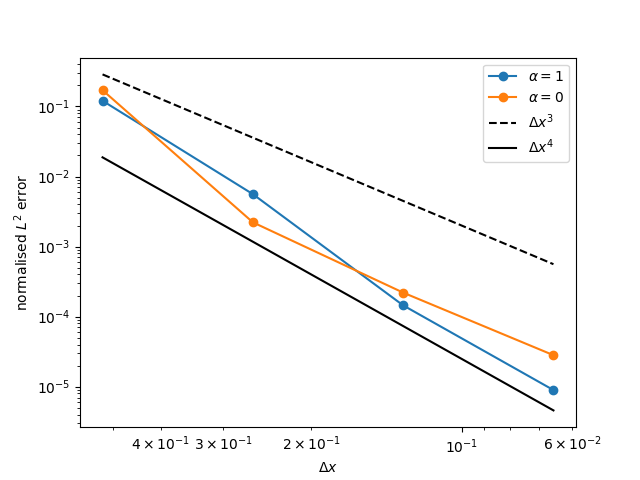}
\includegraphics[width=0.48\textwidth,height=0.36\textwidth]{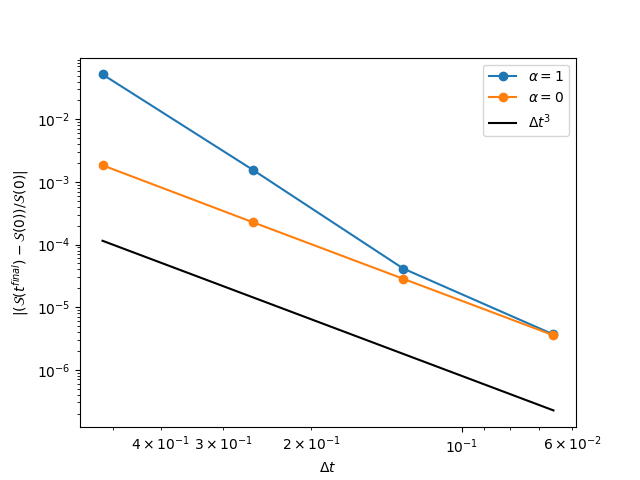}
\caption{$L^2$ error convergence with grid resolution after a single period (left), and absolute tracer variance conservation
	error after a single period (right).}
\label{fig::adv2}
\end{center}
\end{figure}

Since the solid body mass flux is exactly divergence free, mass conservation is assured discretely as per \eqref{eq::adv_mass_con}.
Figure \ref{fig::adv1} gives the tracer mass, $\int h_hs_h\mathrm{d}\Omega$ and tracer variance \eqref{eq::S} 
conservation errors for $h_h=1$, as the normalised difference from their initial values. For the mass 
conservation error there is a small time integration error at each time step such that this decays slightly from machine precision.
For the tracer variance conservation error, this is larger for the upwinded solution ($\alpha=1$) at coarser resolutions, and reduces with increased 
spatial and temporal resolution. 

The convergence of $L^2$ error and tracer variance conservation are given in Fig. \ref{fig::adv2}.
While the $L^2$ errors converge at $4^{\mathrm{th}}$ order for the upwinded flux, for the centered flux there is a decay in the 
convergence rate which is closer to $3^{\mathrm{rd}}$ order.  
Superconvergence of order $p+1$ has been previously established for upwinded variants of discontinuous Galerkin methods for transport 
problems \cite{CockburnEtAl10}. While the stiffly stable Runge-Kutta time 
stepping scheme \cite{SO88} is only $3^{\mathrm{rd}}$ order accurate, the errors incurred from time discretisation are presumably small 
with respect to the spatial discretistion error.
For the tracer variance conservation errors, these decay with the $3^{\mathrm{rd}}$ order stiffly stable Runge-Kutta 
time stepping scheme
as expected for the centered flux, since the spatial discretisation exactly conserves tracer variance in the absence of upwinding.
For the upwinded solution, the absolute value of the tracer variance conservation error is larger as expected for a given time
step size, however the rate of convergence is surprising faster at coarser resolutions.

\subsection{Thermogeostrophic balance}

In order to verify the full low order dynamics/high order transport thermal shallow water solver 
\eqref{eq::tsw_disc_1}, \eqref{eq::tsw_diag_1}, \eqref{eq::buoy_adv_disc},
we perform a convergence 
test for a steady solution to the thermal shallow water equations in thermogeostrophic balance \cite{Eldred19,Tambyah25} using 
periodic boundary conditions. The 
test was run on a domain using comparable scales to that of the earth, with a length of $L=2\pi r_e$ in both dimensions for $r_e=6371220.0$m,
with low order elements of size $\Delta x = L/32,L/64,L/128,L/256$m for a total time of $24$ hours using time steps of $\Delta t=1800, 900, 450, 225$s,
and upwinded high order buoyancy transport ($\alpha=1$).
The initial steady state condition is given as $\boldsymbol{u}=(U_0\cos(y/r_e),0)$, $h=H_0 - r_efU_0/g\sin(y/r_e)$, 
$s=g(1 + 0.05H_0^2/h^2)$, where $U_0=20\mathrm{ms}^{-1}$, $H_0=5960$m, $g=9.80616\mathrm{ms}^{-2}$, $f=6.147\times 10^{-5}\mathrm{s}^{-1}$.

While upwinding is not strictly necessary for this test, owing to the smoothness of the solution, for long times with transient turbulent
flows upwinding may be necessary to maintain both stability and convergence. This is because while the transport scheme conserves tracer variance 
in space, it is not conserved in time, and so in the absence of upwinding this can grow uncontrollably, resulting in grid scale aliasing errors
that may effect convergence or stability.

\begin{figure}[!hbtp]
\begin{center}
\includegraphics[width=0.48\textwidth,height=0.36\textwidth]{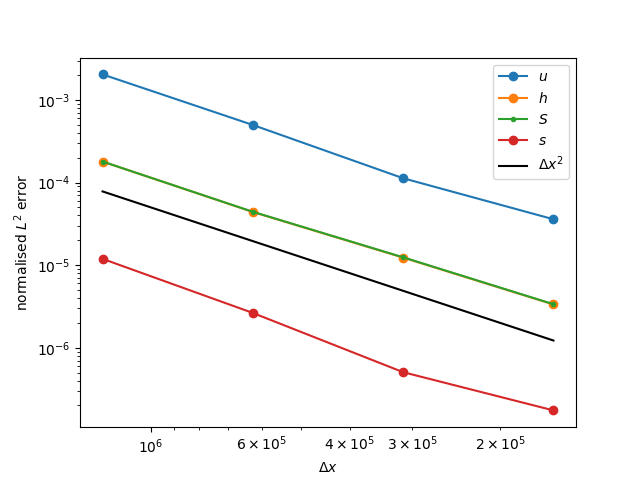}
\includegraphics[width=0.48\textwidth,height=0.36\textwidth]{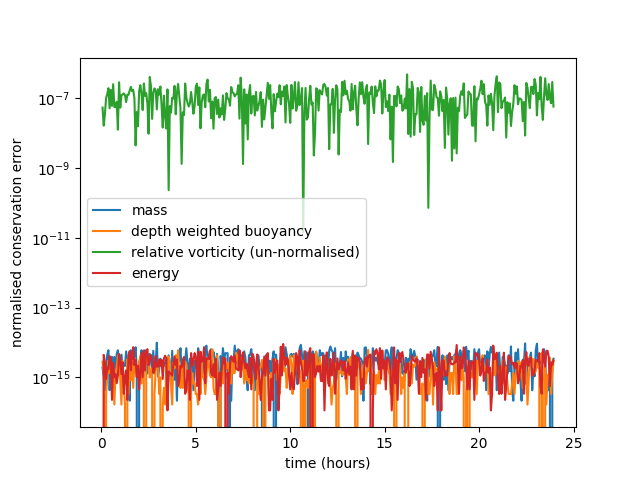}
	\caption{$L^2$ error convergence (left) and normalised conservation errors for the mass, 
	depth weighted buoyancy, relative vorticity and energy using $256\times 256$ low order elements (right).}
\label{fig::thermo_geo1}
\end{center}
\end{figure}

Figure \ref{fig::thermo_geo1} gives the error convergence for the prognostic variables, $\boldsymbol{u}_h, h_h, S_h$,
as well as the buoyancy $s_h$. All variables converge at second order as expected for the low order spatial discretisation with centered time
integration. Unlike for the pure advection test, here the convergence for the buoyancy, $s_h$ is also limited to second order. This is due
to the fact that this is diagnosed from the low order prognostic variables, $h_h, S_h$ at the beginning of each time step, and due to the
low order mass flux used for buoyancy transport. The conservation errors for the mass, depth weighted buoyancy, relative vorticity, 
$\omega=\nabla\times\boldsymbol{u}$ and total energy are also given in Fig. \ref{fig::thermo_geo1} for the highest resolution ($256\times 256$ low order elements). For the mass, depth weighted buoyancy 
and energy, these 
are normalised by their initial value, and are at machine precision. For the relative vorticity, this is un-normalised, 
since the initial value integrates to 0, and are of $\mathcal{O}(10^{-7})$, without exhibiting any long term drift.

\begin{figure}[!hbtp]
\begin{center}
\includegraphics[width=0.48\textwidth,height=0.36\textwidth]{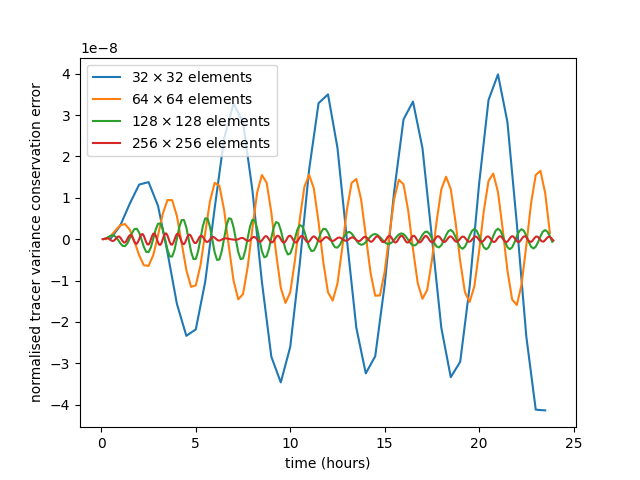}
\includegraphics[width=0.48\textwidth,height=0.36\textwidth]{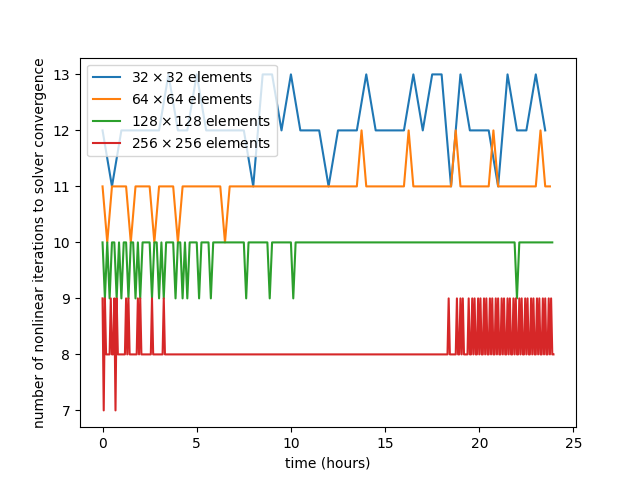}
	\caption{Normalised conservation error with resolution for the tracer variance (left),
	and number of iterations to nonlinear solver convergence with resolution (right).}
\label{fig::thermo_geo2}
\end{center}
\end{figure}

The normalised tracer variance conservation errors are presented for each resolution in Fig. \ref{fig::thermo_geo2}.
These exhibit a small oscillation of $\mathcal{O}(10^{-8})$, which decays in wavelength and amplitude with increasing 
spatial and temporal resolution. 
Figure \ref{fig::thermo_geo2} also shows the number of iterations for the
nonlinear solver as a function of time for each resolution, with convergence given by $||\delta\boldsymbol{u}_h||/||\boldsymbol{u}_h||,
||\delta h_h||/||h_h||, ||\delta S_h||/||S_h|| < 10.0^{-12}$. This reduces for increased spatial and temporal resolution.

\subsection{Thermal instability}

In order to verify the thermal shallow water solver in a more well developed nonlinear regime and to quantify the benefits of the high order
transport, this was applied to a standard test case for a thermally unstable single vortex \cite{Eldred19,KLZ21,Tambyah25}.
The periodic domain was set as $\Omega = [-\pi,\pi)\times[-\pi,\pi)$ and the time step as $\Delta t=0.02$ for a total simulation time
of 100 dimensionless units, using 288 fine scale low order elements (72 $4^{\mathrm{th}}$ order discontinuous Galerkin 
elements) in each dimension.

The initial conditions are given in polar coordinates $r=\sqrt{x^2+y^2},\theta=\tan^{-1}(y/x)$ as a perturbation, $\epsilon$
to a balanced state of 
\begin{subequations}
	\begin{align}
		u&=\epsilon - U_0r\exp((1-r^{\beta})/\beta)\sin(\theta),\\
		v&=\epsilon + U_0r\exp((1-r^{\beta})/\beta)\cos(\theta),\\
		h&=H_0-\epsilon,\\
		s&=\epsilon+g-\frac{2R_o}{B_u}\Big(\exp((1-r^2)/2) + \frac{R_o}{2}\exp(1-r^2)\Big),
	\end{align}
\end{subequations}
with the perturbation as
\begin{equation}
	\epsilon=0.01\exp(-60(r-r_c)^2)\sin(6\pi(r-r_c))\cos(4\theta).
\end{equation}
The constants are given as $H_0=g=1$, $\beta=2$, $r_c=0.5$,
$U_0=0.1$ and the Rossby and Burgers numbers respectively as $R_o=0.1$, $B_u=1$.
For this test the solver was configured to run for a fixed number of four Newton iterations per time step, rather than to 
convergence, such that energy conservation was not preserved in time, in order to better reflect the configuration in
operational atmospheric models, where it is not efficient to run the nonlinear solver to convergence. 

The high order discontinuous Galerkin buoyancy configuration with upwinding ($\alpha=1$) is compared to three alternative low order formulations as follows,
none of which use the high order transport scheme:
\begin{itemize}
	\item \emph{Low order centered flux form transport:} 
		Where the buoyancy is derived solely from the low order dynamics subject to the low order test function
		$\phi_h\in\mathbb{W}_2^L$ as 
		\begin{equation}
		\int\phi_hh_hs_h\mathrm{d}\Omega=\int\phi_hS_h\mathrm{d}\Omega.
		\end{equation}
	\item \emph{Low order centered skew-symmetric flux form transport:}
		The buoyancy is derived as above, however the pressure gradient and buoyancy flux terms 
		are re-formulated within \eqref{eq::tsw_disc::u} and \eqref{eq::tsw_disc::S} respectively as the
		low order analogue of the tracer and energy conserving formulation \cite{Ricardo24a,Tambyah25} as
\begin{subequations}\label{eq::tsw_disc_lo}
\begin{align}
		\int\boldsymbol{v}_h\cdot(\boldsymbol{u}_h^{n+1} - \boldsymbol{u}_h^n)\mathrm{d}\Omega + 
		\Delta t\int\boldsymbol{v}_h\cdot\overline{q}_h\overline{\boldsymbol{F}}_h^{\perp}\mathrm{d}\Omega 
		-\Delta t\int\nabla\cdot\boldsymbol{v}_h\overline{\Phi}_h\mathrm{d}\Omega \notag\\ 
		-\frac{\Delta t}{2}\int\boldsymbol{v}_h\cdot\hat{\boldsymbol{n}}\{\overline{s}_h\}[\overline{T}_h]\mathrm{d}\Gamma 
		+\frac{\Delta t}{2}\int\boldsymbol{v}_h\cdot\hat{\boldsymbol{n}}[\overline{s}_h]\{\overline{T}_h\}\mathrm{d}\Gamma 
		-\frac{\Delta t}{2}\int\nabla\cdot\boldsymbol{v}_h\overline{s}_h\overline{T}_h\mathrm{d}\Omega &= 0,\label{eq::disc_u_low_tvc}\\
		\int\sigma_h(S_h^{n+1} - S_h^n)\mathrm{d}\Omega 
		+\frac{\Delta t}{2}\int[\sigma_h]\{\overline{s}_h\}\overline{\boldsymbol{F}}_h\cdot\hat{\boldsymbol{n}}\mathrm{d}\Gamma 
		-\frac{\Delta t}{2}\int\{\sigma_h\}[\overline{s}_h]\overline{\boldsymbol{F}}_h\cdot\hat{\boldsymbol{n}}\mathrm{d}\Gamma 
		+\frac{\Delta t}{2}\int\sigma_h\overline{s}_h\nabla\cdot\overline{\boldsymbol{F}}_h\mathrm{d}\Omega
		&= 0.\label{eq::tsw_disc_lo::S}
\end{align}
\end{subequations}
	\item \emph{Low order upwinded skew-symmetric flux form transport:}
This configuration is the same as that above \eqref{eq::tsw_disc_lo}, with the addition of an 
upwinding term to \eqref{eq::tsw_disc_lo::S} of the form
$\Delta t\int|\overline{\boldsymbol{F}}_h\cdot\hat{\boldsymbol{n}}|[\sigma_h][\overline{s}_h]\mathrm{d}\Gamma$
so as to explicitly dissipate tracer variance as a low order analogue to \eqref{eq::disc_var_con}.
One could also add the adjoint of this term to \eqref{eq::disc_u_low_tvc} in order to preserve energy conservation \cite{Eldred19}.
\end{itemize}

\begin{figure}[!hbtp]
\begin{center}
\includegraphics[width=0.48\textwidth,height=0.36\textwidth]{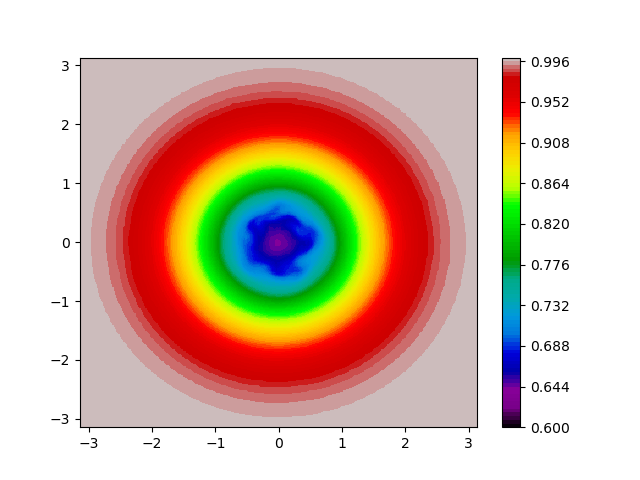}
\includegraphics[width=0.48\textwidth,height=0.36\textwidth]{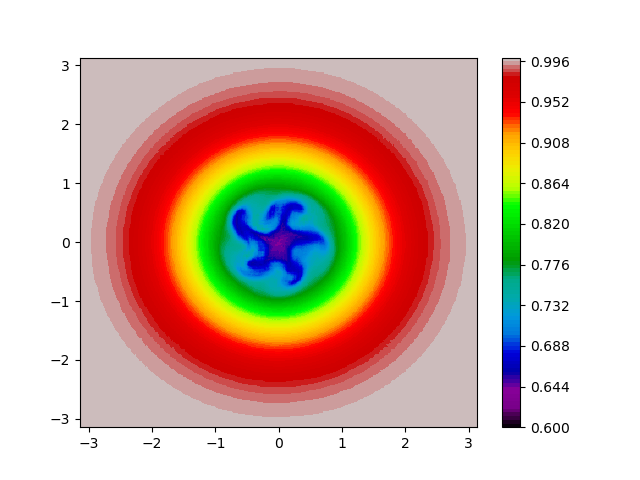}
\includegraphics[width=0.48\textwidth,height=0.36\textwidth]{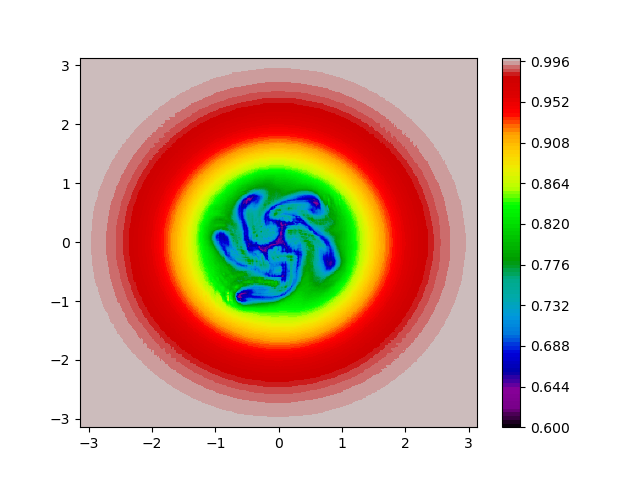}
\includegraphics[width=0.48\textwidth,height=0.36\textwidth]{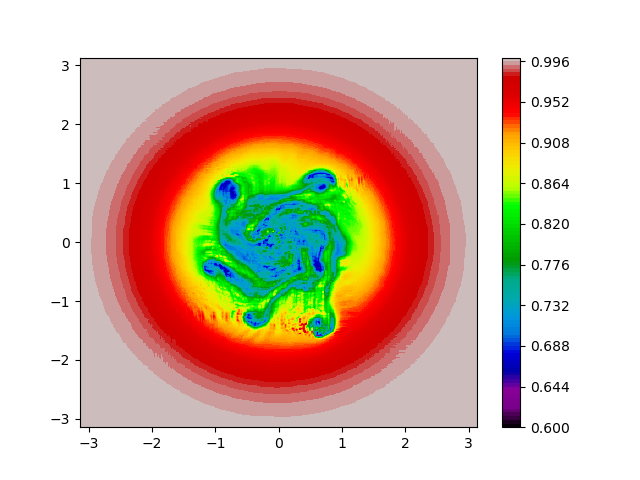}
	\caption{Buoyancy field for the thermal instability test case using high order upwinded transport 
	at dimensionless times $20$ (top left), $40$ (top right), $60$ (bottom left), $80$ (bottom right).}
\label{fig::thermal_instability_1}
\end{center}
\end{figure}

\begin{figure}[!hbtp]
\begin{center}
\includegraphics[width=0.48\textwidth,height=0.36\textwidth]{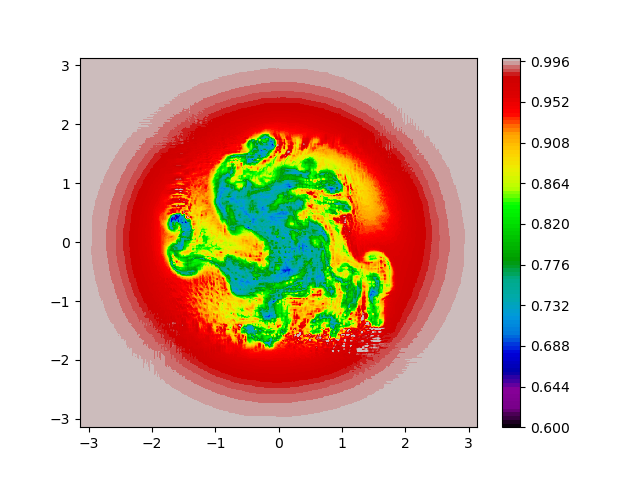}
\includegraphics[width=0.48\textwidth,height=0.36\textwidth]{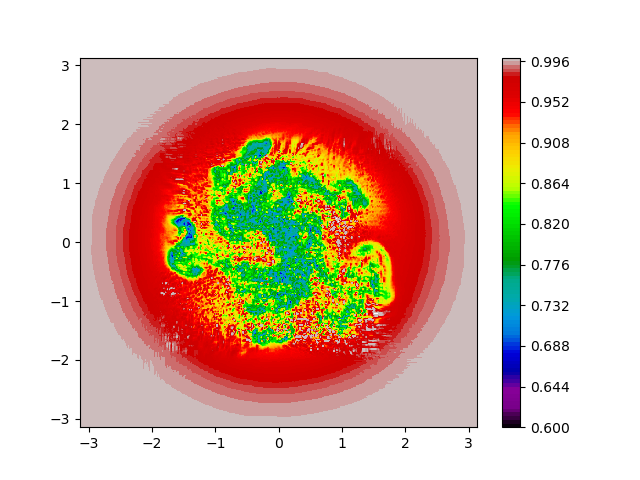}
\includegraphics[width=0.48\textwidth,height=0.36\textwidth]{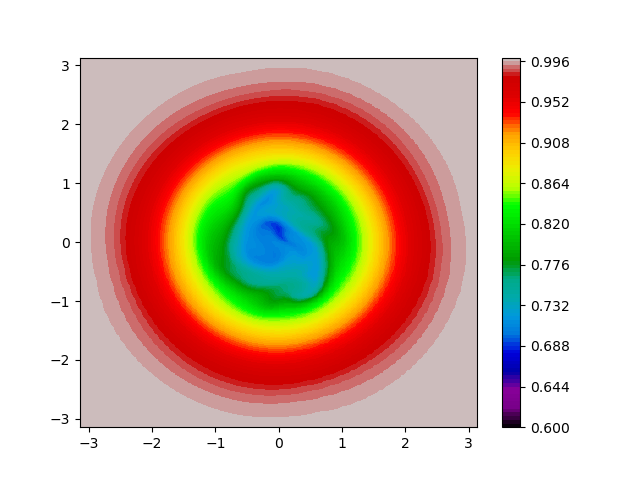}
\includegraphics[width=0.48\textwidth,height=0.36\textwidth]{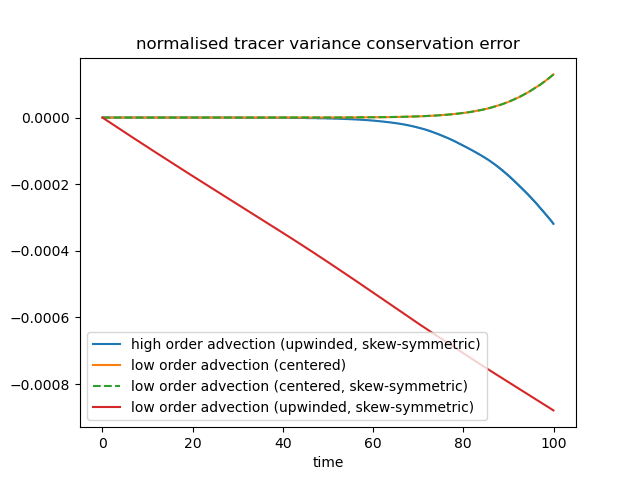}
	\caption{Buoyancy field for the thermal instability test case at dimensionless time $100.0$
	using high order upwinded buoyancy advection (top left), low order skew-symmetric centered advection (top right), 
	low order skew-symmetric upwinded advection (bottom left), and normalised tracer variance
	conservation errors for the thermal instability test case (bottom right).}
\label{fig::thermal_instability_2}
\end{center}
\end{figure}

The buoyancy field for the original configuration \eqref{eq::tsw_disc_1}, \eqref{eq::tsw_diag_1}, \eqref{eq::buoy_adv_disc} 
using high order discontinuous Galerkin transport of $s_h$ is given
at dimensionless times of $20,40,60,80$ in Fig. \ref{fig::thermal_instability_1}. The corresponding solution at time
$100$ is given in comparison to those using low order skew-symmetric buoyancy transport (both centered and upwinded) in Fig. 
\ref{fig::thermal_instability_2}. The low order upwinded solution is excessively diffusive, which is expected, 
since the inclusion of the low order upwinding term effectively means that the low order flux is sampled only 
from the (piecewise constant) upwind cell, meaning that this flux is effectively first order only. The low order centered
flux by contrast is excessively noisy, since while the tracer variance is conserved (in space but not in time) for this
formulation, there is nothing to control the nonlinear aliasing of the tracer variance at the grid scale.

These results are also reflected in the tracer variance conservation error time series for the three different configurations,
which are also given in Fig. \ref{fig::thermal_instability_2}. While the tracer variance conservation error for the low order 
upwinded solution decays linearly, reflecting the excessively diffusive solution, for the low order centered flux this 
grows with time, suggesting that this solution will eventually become unstable. 
This is despite the fact that tracer variance
is conserved by the spatial (but not the temporal) discretisation. By contrast, the tracer variance conservation
error for the high order discontinuous Galerkin buoyancy transport (which also includes upwinding) decays, suggesting model 
stability, but at a more moderate rate than for the low order upwinded solution.
The results are identical (to machine precision) for the skew-symmetric
and non-skew-symmetric centered flux formulations, which demonstrates that for the lowest order formulation, this is indeed 
tracer variance conserving in space without the skew-symmetric correction as discussed in Section 2.3 \eqref{eq::lo_tv_disc}.

\section{Conclusions}

This article presents a coupling of a low order, mixed finite element formulation of the thermal shallow
water equations with a high order, tracer variance conserving discontinuous Galerkin scheme for buoyancy advection in
material form.
The degrees of freedom for the high order elements on which the buoyancy transport is computed are collocated 
with the centres of the low order elements for the dynamics at the Gauss Legendre points on the high order mesh.
The high order tracer variance conserving transport 
scheme extends previous work on the formulation of tracer variance conserving methods for flux form advection
using discontinuous Galerkin \cite{Ricardo24b,Ricardo25} and mixed finite element methods \cite{Ricardo24a,Tambyah25}
by presenting an analogous formulation for material transport. 

Energy conservation is preserved for the low order solver since the high order buoyancy is applied so as to preserve
the antisymmetric structure of the Hamiltonian form of the equations of motion. While the overall method is limited 
to second order accuracy due to the diagnosis of the buoyancy from the low order dynamics at each time step, numerical 
experiments for well developed turbulence in a thermal instability test case show that the use of high order buoyancy
transport with low order dynamics ensures that the solution is nonlinearly stable with respect to tracer variance without
being excessively diffusive, as is the case for upwinding of the low order solution.

While this coupling strategy is presented here for the thermal shallow water equations, we note that the same idea can be applied
to conserve energy and conserve or bound tracer variance for other non-canonical Hamiltonian systems involving the material transport 
of thermodynamic tracers, such as potential temperature or thermodynamic entropy and moisture fractions in the case of the dry \cite{LP21} 
and moist \cite{Ricardo25} compressible Euler equations.

The scheme may be further extended by incorporating the high order buoyancy transport mesh within a hierarchy of meshes as
part of a geometric multigrid solver for the low order dynamics, which may exhibit improved scalability on parallel machines
when used in conjunction with an appropriate patch smoother compared to the existing GMRES Krylov solver used in this study.

\section{Acknowledgements}

David Lee would like to thank Dr. Junwei Lyu for his careful and constructive review of an early version of this manuscript.
Kieran Ricardo would like to acknowledge the Australian Government through the Australian Government Research Training Program (RTP) Scholarship, and the Bureau of Meteorology through research contract KR2326.
Tamara Tambyah would like to acknowledge the Commonwealth of Australia as represented by the Defence Science and Technology Group of the Department of Defence (agreement number 11652).

\end{document}